\documentclass[fleqn,preprint,3p,a4paper]{elsarticle}

\usepackage{amssymb,mathrsfs,latexsym,amsfonts}
\usepackage{amsmath}
\usepackage{amsthm}
\usepackage{dcolumn}
\usepackage{endnotes}
\usepackage{tabularx}
\usepackage[matrix,arrow]{xy}
\usepackage{wasysym}

\newtheorem{theorem}{Theorem}[section]
\newtheorem{proposition}[theorem]{Proposition}
\newtheorem{lemma}[theorem]{Lemma}
\newtheorem{corollary}[theorem]{Corollary}

\theoremstyle{definition}
\newtheorem{definition}[theorem]{Definition}
\newtheorem{example}[theorem]{Example}
\newtheorem{remark}[theorem]{Remark}
\newtheorem{question}[theorem]{Question}

\journal{Journal}

\begin{document}

\begin{frontmatter}

%% Title, authors and addresses

%% use the tnoteref command within \title for footnotes;
%% use the tnotetext command for theassociated footnote;
%% use the fnref command within \author or \address for footnotes;
%% use the fntext command for theassociated footnote;
%% use the corref command within \author for corresponding author footnotes;
%% use the cortext command for theassociated footnote;
%% use the ead command for the email address,
%% and the form \ead[url] for the home page:
%% \title{Title\tnoteref{label1}}
%% \tnotetext[label1]{}
%% \author{Name\corref{cor1}\fnref{label2}}
%% \ead{email address}
%% \ead[url]{home page}
%% \fntext[label2]{}
%% \cortext[cor1]{}
%% \address{Address\fnref{label3}}
%% \fntext[label3]{}

\title{On some kinds of factorizable topological groups\tnoteref{t1}}
\tnotetext[t1]{This research was supported by the National Natural Science Foundation of China (Nos. 12071199, 11661057), the Natural Science Foundation of Jiangxi Province, China (No. 20192ACBL20045).}

%% use optional labels to link authors explicitly to addresses:
%% \author[label1,label2]{}
%% \address[label1]{}
%% \address[label2]{}

\author[M. Bao]{Meng Bao}
\ead{mengbao95213@163.com}
\address[M. Bao]{College of Mathematics, Sichuan University, Chengdu 610064, China}

\author[X. Xu]{Xiaoquan Xu\corref{mycorrespondingauthor}}
\cortext[mycorrespondingauthor]{Corresponding author.}
\ead{xiqxu2002@163.com}
\address[X. Xu]{Fujian Key Laboratory of Granular Computing and Applications, Minnan Normal University, Zhangzhou 363000, China}

\begin{abstract}
Based on the concepts of $\mathbb{R}$-factorizable topological groups and $\mathcal{M}$-factorizable topological groups, we introduce four classes of factorizabilities on topological groups, named $P\mathcal{M}$-factorizabilities, $Pm$-factorizabilities, $S\mathcal{M}$-factorizabilities and $PS\mathcal{M}$-factorizabilities, respectively. Some properties of the four classes of spaces are investigated.
\end{abstract}


%%Graphical abstract
%\begin{graphicalabstract}
%\includegraphics{grabs}
%\end{graphicalabstract}

%%Research highlights
%\begin{highlights}
%\item Research highlight 1
%\item Research highlight 2
%\end{highlights}

%% MSC codes here, in the form:
%% or \MSC[2008] code \sep code (2000 is the default)

\end{frontmatter}

%% \linenumbers

%% main text
\section{Introduction}

In the field of Topological Algebra, topological groups are standard researching objects and have been studied for many years, see \cite{AA}. A {\it topological group} is a group equipped with a topology such that the multiplication on the group is jointly continuous and the inverse mapping is also continuous. It is well-known that for every continuous real-valued function $f$ on a compact topological group $G$, there exists a continuous homomorphism $p:G\rightarrow L$ onto a second-countable topological group $L$ and a continuous real-valued function $h$ on $L$ such that $f=h\circ p$. Then, Tkachenko posed the concept of $\mathbb{R}$-factorizable topological groups, see \cite{Tkachenko91}. A topological group $G$ is called {\it $\mathbb{R}$-factorizable} if for every continuous real-valued function $f$ on $G$, we can find a continuous homomorphism $\pi :G\rightarrow H$ onto a second-countable topological group $H$ such that $f=g\circ \pi$, for some continuous real-valued function $g$ on $H$. We know that $\mathbb{R}$-factorizable topological groups are generalizations of compact groups and separable metrizable groups. For more interesting properties about $\mathbb{R}$-factorizable topological groups, see \cite{Tkachenko2006, Tkachenko2010, Tkachenko2012}. However, since a metrizable topological group need not to be $\mathbb{R}$-factorizable, it follows that H. Zhang, D. Peng and W. He in \cite{ZH2020} posed the notion of $\mathcal{M}$-factorazible topological group. A topological group $G$ is called {\it $\mathcal{M}$-factorizable} if for every continuous real-valued function $f$ on $G$, there is a continuous homomorphism $\varphi :G\rightarrow H$ onto a metrizable topological group $H$ such that $f=g\circ \varphi$, for some continuous real-valued function $g$ on $H$. Since all first-countable topological groups are metrizable, it is trivial to see that all first-countable topological groups are $\mathcal{M}$-factorazible. Moreover, it was proved in \cite[Theorem 3.2]{ZH2020} that a topological group is $\mathbb{R}$-factorizable if and only if it is $\mathcal{M}$-factorazible and $\omega$-narrow.

By further researches of $\mathbb{R}$-factorizable topological groups, L. Peng and Y. Liu introduced the concept of $P\mathbb{R}$-factorizable topological groups, that is, a topological group $G$ is called {\it $P\mathbb{R}$-factorizable} if for every continuous real-valued function $f$ on $G$, there exists a perfect homomorphism $\pi :G\rightarrow H$ onto a second-countable topological group $H$ such that $f=g\circ \pi$, for some continuous real-valued function $g$ on $H$. They gave the characterizations of $P\mathbb{R}$-factorizable topological groups in \cite[Theorem 2.6]{Peng2020}. In particular, a topological group is $P\mathbb{R}$-factorizable if and only if it is Lindel\"{o}f feathered. Moreover, since every $\omega$-narrow feathered topological group is Lindel\"{o}f (see \cite[4.3.A]{AA}), it is easy to see that a topological group is $P\mathbb{R}$-factorizable if{}f it is $\omega$-narrow and feathered. Then, we introduce the following notion.

\begin{definition}
A topological group $G$ is called {\it $P\mathcal{M}$-factorizable} if for every continuous real-valued function $f$ on $G$, there exists a perfect homomorphism $\pi :G\rightarrow H$ onto a metrizable topological group $H$ such that $f=g\circ \pi$, for some continuous real-valued function $g$ on $H$.
\end{definition}

Clearly, each $P\mathbb{R}$-factorizable topological group is $P\mathcal{M}$-factorizable. L. Peng and Y. Liu introduced an example \cite[Example 3.13]{Peng2020} which is an $\mathbb{R}$-factorizable topological group, but not $P\mathbb{R}$-factorizable. Indeed, the topological group $G$ in \cite[Example 3.13]{Peng2020} is not feathered. Since all $\mathbb{R}$-factorizable topological groups are $\mathcal{M}$-factorizable, it is a $\mathcal{M}$-factorizable topological group. However, by the definition of $\mathcal{M}$-factorizability, it is easy to that every $P\mathcal{M}$-factorizable topological group is feathered, hence the topological group $G$ of \cite[Example 3.13]{Peng2020} is not $P\mathcal{M}$-factorizable.

In this paper, we give some characterizations of $P\mathcal{M}$-factorizable topological groups, such as a topological group $G$ is $P\mathcal{M}$-factorizable if and only if $G$ is feathered $\mathcal{M}$-factorizable. We also shown that a topological group $G$ is $P\mathbb{R}$-factorizable if and only if $G$ is $P\mathcal{M}$-factorizable and $\omega$-narrow. Then it is natural to deduce that a topological group $G$ is $P\mathbb{R}$-factorizable if and only if $G$ is feathered $\mathbb{R}$-factorizable. W. He et al. in \cite[Proposition 2.1]{HW2021} proved the following proposition.

\begin{proposition}
Let $G=\prod_{i\in I}G_{i}$ be the product of an uncountable family of non-compact separable metrizable topological groups. Then the group $G$ is $\mathbb{R}$-factorizable, but it fails to be feathered
\end{proposition}

Therefore, the result also can present that the product group $G$ is not $P\mathbb{R}$-factorizable. Of course, $G$ is an $\mathcal{M}$-factorizable topological group, but not $P\mathcal{M}$-factorizable. Moreover, some interesting properties of $\mathcal{M}$-factorizable topological groups in \cite{HW2021, ZH2020} are strengthened to $P\mathcal{M}$-factorizable topological groups. For example, the product $G=\prod_{n\in \mathbb{N}}G_{n}$ of countably many $P\mathcal{M}$-factorizable topological groups is $P\mathcal{M}$-factorizable if and only if every factor $G_{n}$ is metrizable or every $G_{n}$ is $P\mathbb{R}$-factorizable, the product of a $P\mathcal{M}$-factorizable topological group with a compact metrizable topological group is $P\mathcal{M}$-factorizable.

Then, according to the concept of $m$-factorizable topological groups, that is, a topological group $G$ is called {\it $m$-factorizable} if for every continuous mapping $f:G\rightarrow M$ to a metrizable space $M$, there exists a continuous homomorphism $\pi :G\rightarrow K$ onto a second-countable group $K$ such that $f=g\circ \pi$, for some continuous homomorphism $g$ from $K$ onto $M$, we introduce $Pm$-factorizable topological groups by strengthening the continuous homomorphism $\pi$ to a perfect homomorphism. Then we show that a topological group $G$ is $Pm$-factorizable if and only if $G$ is $P\mathbb{R}$-factorizable and pseudo-$\aleph_{1}$-compact, which deduces that a topological group $G$ is $Pm$-factorizable if and only if $G$ is feathered $m$-factorizable. Furthermore, it is claimed that the product of a $Pm$-factorizable topological group with an arbitrary compact group is $Pm$-factorizable.

Finally, we pose the concepts of $S\mathcal{M}$-factorizable topological groups and $PS\mathcal{M}$-factorizable topological groups. A topological group $G$ is called {\it strongly $\mathcal{M}$-factorizable} ($S\mathcal{M}$-factorizable for short) if for every continuous mapping $f:G\rightarrow M$ to a metrizable space $M$, there exists a continuous homomorphism $\pi :G\rightarrow H$ onto a metrizable group $H$ and a continuous mapping $g:H\rightarrow M$ such that $f=g\circ \pi$. In particular, if the continuous homomorphism $\pi :G\rightarrow H$ is perfect, we call $G$ {\it $PS\mathcal{M}$-factorizable.} Since a topological group $G$ is $P\mathcal{M}$-factorizable if and only if $G$ is feathered $\mathcal{M}$-factorizable, and a topological group $G$ is $Pm$-factorizable if and only if $G$ is feathered and $m$-factorizable, see Theorem \ref{PM=M} and Corollary \ref{Pm=m}, it is natural to consider whether is it equivalent between $PS\mathcal{M}$-factorizability and feathered $S\mathcal{M}$-factorizability. Then we show that it also holds in Theorem \ref{PSM=feathered+SM}. Furthermore, the four classes of factorizable properties are preserved by open continuous homomorphisms on topological groups, see Theorem \ref{PMopen+continuous}, Corollary \ref{Pmopen+continuous}, Theorem \ref{SMopen+continuous} and Corollary \ref{PSMopen+continuous}.

\section{Preliminary}

Throughout this paper, all topological spaces are assumed to be Hausdorff, unless otherwise is explicitly stated. Let $\mathbb{N}$ be the set of all positive integers and $\omega$ the first infinite ordinal. The readers may consult \cite{AA, E} for notation and terminology not explicitly given here. Next we recall some definitions and facts.

A continuous mapping $f:X\rightarrow Y$ is called {\it perfect} if it is a closed onto mapping and all fibers $f^{-1}(y)$ are compact subsets of $X$ \cite[p. 182]{E}. A Tychonoff topological space $X$ is {\it \v{C}ech-complete} if $X$ is a $G_{\delta}$-set in every compactification $cX$ of the space $X$ \cite[p. 192]{E}.

Then we recall some notions about topological groups. A topological group $G$ is called {\it feathered} if it contains a non-empty compact set $K$ of countable character in $G$. It is well-known in \cite[Theorem 4.3.20]{AA} that a topological group $G$ is feathered if and only if it contains a compact subgroup $H$ such that the left quotient space $G/H$ is metrizable. Similarly, the group $G$ is \v{C}ech-complete if and only if it contains a compact subgroup $H$ such that the left quotient space $G/H$ is metrizable by a complete metric. The class of feathered topological group is countably productive and is closed under closed-heredity. Moreover, all \v{C}ech-complete topological groups are feathered and every metrizable (or locally compact) topological group is feathered.

A topological group is {\it Ra\v\i kov complete} if it is complete with respect to its two-sided uniform group structure. Every topological group $G$ can be embedded into a unique Ra\v\i kov complete topological group as a dense subgroup, which is called the {\it Ra\v\i kov completion} of $G$ and denoted by $\varrho G$. Topological groups with compact completions are called {\it precompact}. As we all know, every \v{C}ech-complete topological group is Ra\v\i kov complete and a feathered topological group is Ra\v\i kov complete if and only if it is \v{C}ech-complete, see \cite[Theorem 4.3.15]{AA}.

We call a topological group $G$ {\it $\omega$-narrow } if for any neighborhood $U$ of the identity $e$ in $G$, there exists a countable subset $C$ of $G$ such that $G=UC=CU$. A topological group $G$ is called {\it $\omega$-balanced} if for each neighborhood $U$ of the identity $e$ in $G$, there exists a countable family $\gamma$ of open neighborhoods of $e$ such that for each $x\in G$, there exists $V\in \gamma$ satisfying $xVx^{-1}\subseteq U$. \cite[Proposition 3.4.10]{AA} showed that every $\omega$-narrow topological group is $\omega$-balanced.

\vspace{-0.2cm}

\section{Some properties of $P\mathcal{M}$-factorizable topological groups}

In this section, we give some characterizations of $P\mathcal{M}$-factorizable topological groups, such as a topological group $G$ is $P\mathcal{M}$-factorizable if and only if $G$ is feathered $\mathcal{M}$-factorizable. We also shown that a topological group $G$ is $P\mathbb{R}$-factorizable if and only if $G$ is $P\mathcal{M}$-factorizable and $\omega$-narrow. Then it is natural to deduce that a topological group $G$ is $P\mathbb{R}$-factorizable if and only if $G$ is feathered $\mathbb{R}$-factorizable.

Recall that {\it paracompact $p$-spaces} are the preimages of metrizable spaces under perfect mappings. By the definitions of $P\mathcal{M}$-factorizable topological groups, the following result is clear.

\begin{proposition}
Every $P\mathcal{M}$-factorizable topological group is a paracompact $p$-space.
\end{proposition}

It was proved in \cite[Theorem 4.3.35]{AA} that a topological group is feathered if{}f it is a $p$-space, and if{}f it is a paracompact $p$-space.

\begin{proposition}\label{PM-feathered}
Every $P\mathcal{M}$-factorizable topological group is feathered.
\end{proposition}

Then according to the concept of $P\mathbb{R}$-factorizable topological groups, every $P\mathbb{R}$-factorizable topological group is $P\mathcal{M}$-factorizable. Therefore, every compact topological group is $P\mathcal{M}$-factorizable. Moreover, it is not difficult to see that a feathered $\mathcal{M}$-factorizable topological group is $P\mathcal{M}$-factorizable.

\begin{theorem}\label{PM=M}
A topological group $G$ is $P\mathcal{M}$-factorizable if and only if $G$ is feathered $\mathcal{M}$-factorizable.
\end{theorem}

\begin{proof}
The necessity is trivial.

The sufficiency of Theorem \ref{PM=M} can be shown as follows. If $G$ is a feathered $\mathcal{M}$-factorizable topological group, it follows from \cite[Theorem 3.3]{HW2021} that either $G$ is metrizable or $G$ is $\mathbb{R}$-factorizable. If $G$ is metrizable, it is trivial. On the other hand, if $G$ is an $\mathbb{R}$-factorizable topological group, $G$ is $\omega$-narrow by \cite[Proposition 8.1.3]{AA}. Since a topological group is $P\mathbb{R}$-factorizable if and only if it is $\omega$-narrow and feathered by \cite[Theorem 2.6]{Peng2020}. Therefore, $G$ is a $P\mathbb{R}$-factorizable topological group, and is also $P\mathcal{M}$-factorizable, naturally.
\end{proof}

Moreover, we show that each $\omega$-narrow $P\mathcal{M}$-factorizable topological group is $P\mathbb{R}$-factorizable.

\begin{theorem}\label{omega-PM}
A topological group $G$ is $P\mathbb{R}$-factorizable if and only if $G$ is $P\mathcal{M}$-factorizable and $\omega$-narrow.
\end{theorem}

\begin{proof}
Since every $P\mathbb{R}$-factorizable topological group is $\omega$-narrow, the necessity is trivial.

Let's prove the sufficiency. Suppose that $G$ is a $P\mathcal{M}$-factorizable and $\omega$-narrow topological group, $f:G\rightarrow \mathbb{R}$ is a continuous real-valued function. Then there exists a perfect homomorphism $\varphi :G\rightarrow K$ onto a metrizable topological group $K$ such that $f=g\circ \varphi$, where $g:K\rightarrow \mathbb{R}$ is continuous. Since $G$ is $\omega$-narrow, so is $K$. Since every first-countable $\omega$-narrow topological group is second-countable, we obtain that $G$ is $P\mathbb{R}$-factorizable.
\end{proof}

By Theorems \ref{PM=M} and \ref{omega-PM}, we obtain the following result.

\begin{corollary}\label{PR=R}
A topological group $G$ is $P\mathbb{R}$-factorizable if and only if $G$ is feathered $\mathbb{R}$-factorizable.
\end{corollary}

Indeed, it was proved in \cite[Theorem 3.4]{HW2021} that for a feathered topological group $G$, $G$ is $\omega$-narrow if and only if it is $\mathbb{R}$-factorizable. Moreover, \cite[Theorem 2.6]{Peng2020} presented that a topological group $G$ is $P\mathbb{R}$-factorizable if and only if $G$ is $\omega$-narrow and feathered. As a topological group $G$ is $\mathbb{R}$-factorizable if and only if it is $\mathcal{M}$-factorizable and $\omega$-narrow by \cite[Theorem 3.2]{ZH2020}, the Corollary \ref{PR=R} also can be obtained.

\begin{proposition}
A topological group $G$ is $P\mathbb{R}$-factorizable if and only if $G$ is a feathered Lindel\"{o}f $\Sigma$-group.
\end{proposition}

\begin{proof}
The sufficiency is clear. Indeed, a feathered Lindel\"{o}f $\Sigma$-group $G$ is Lindel\"{o}f feathered, so $G$ is $P\mathbb{R}$-factorizable, as a topological group is Lindel\"{o}f feathered if and only if it is $P\mathbb{R}$-factorizable by \cite[Theorem 2.5]{Peng2020}.

Then we show the necessity. Let $G$ be a $P\mathbb{R}$-factorizable topological group. It follows from \cite[Theorem 2.6]{Peng2020} that $G$ is $\omega$-narrow and feathered. Then, by \cite[Theorem 3.4]{HW2021}, for a feathered topological group $G$, $G$ is $\omega$-narrow if and only if $G$ is a Lindel\"{o}f $\Sigma$-group.
\end{proof}

\begin{theorem}\label{metrizable-PR}
A topological group $G$ is $P\mathcal{M}$-factorizable if and only if one of the following holds:
\begin{enumerate}
\smallskip
\item $G$ is metrizable;

\smallskip
\item $G$ is $P\mathbb{R}$-factorizable.
\end{enumerate}
\end{theorem}

\begin{proof}
Since every $P\mathbb{R}$-factorizable topological group is $P\mathcal{M}$-factorizable and every metrizable topological group is also $P\mathcal{M}$-factorizable, the sufficiency is clear.

Then suppose that a $P\mathcal{M}$-factorizable topological group $G$ is not metrizable. By Proposition \ref{PM-feathered}, $G$ is feathered. Since a non-metrizable $P\mathcal{M}$-factorizable topological group is $\mathbb{R}$-factorizable by \cite[Theorem 3.3]{HW2021}, $G$ is $\mathbb{R}$-factorizable and so $G$ is a $P\mathbb{R}$-factorizable topological group
\end{proof}

Then from \cite[Theorem 4.8]{ZH2020}, a locally compact group $G$ is $\mathcal{M}$-factorizable if and only if $G$ is metrizable or $G$ is $\sigma$-compact. Then, it is well-known that every locally compact topological group is feathered, so the characterization also holds for $P\mathcal{M}$-factorizable topological groups by Theorem \ref{PM=M}.

\begin{proposition}\label{3mt10}
A locally compact group $G$ is $P\mathcal{M}$-factorizable if and only if one of the following conditions holds:
\begin{enumerate}
\smallskip
\item $G$ is metrizable;

\smallskip
\item $G$ is $\sigma$-compact.
\end{enumerate}
\end{proposition}

It is well-known that a subgroup of an $\mathbb{R}$-factorizable topological group may not be $\mathbb{R}$-factorizable. Indeed, by \cite[Example 8.2.1]{AA}, there is an Abelian $P$-group $G$ and a dense subgroup $H$ of $G$ such that $G$ is Lindel\"{o}f, hence $\mathbb{R}$-factorizable, but $H$ is not $\mathbb{R}$-factorizable. In particular, $H$ is $\omega$-narrow. Therefore, $H$ is not $\mathcal{M}$-factorizable. Then since every $\omega$-narrow topological group can be embedded into an $\mathbb{R}$-factorizable as a closed invariant subgroup, see \cite[Theorem 8.2.2]{AA}, hence $\mathbb{R}$-factorizability is not closed-heredity for topological groups. However, W. He et al. showed that every subgroup of an $\mathcal{M}$-factorizable feathered group is $\mathcal{M}$-factorizable, it also means that every subgroup of a $P\mathcal{M}$-factorizable is $\mathcal{M}$-factorizable. So it is easy to see that $P\mathcal{M}$-factorizability is closed-heredity for topological groups.

\begin{proposition}
Every closed subgroup of a $P\mathcal{M}$-factorizable topological group is $P\mathcal{M}$-factorizable.
\end{proposition}

\begin{proof}
Let $G$ be a $P\mathcal{M}$-factorizable topological group and $H$ a closed subgroup of $G$. By Theorem \ref{PM=M}, $G$ is $\mathcal{M}$-factorizable and feathered. According to \cite[Lemma 4.1]{HW2021}, every subgroup of an $\mathcal{M}$-factorizable feathered group is $\mathcal{M}$-factorizable, so $H$ is a $\mathcal{M}$-factorizable topological group. Moreover, it is well-known that a closed subspace of a feathered space is feathered. Hence, $H$ is $\mathcal{M}$-factorizable and feathered. We have that $H$ is $P\mathcal{M}$-factorizable by Theorem \ref{PM=M}.
\end{proof}

From \cite[Theorem 3.4.4]{AA}, every subgroup of an $\omega$-narrow topological group is $\omega$-narrow, so the following is clear.

\begin{corollary}\cite[Proposition 2.8]{Peng2020}
Every closed subgroup of a $P\mathbb{R}$-factorizable topological group is $P\mathbb{R}$-factorizable.
\end{corollary}

\begin{theorem}\label{PMopen+continuous}
If $f:G\rightarrow H$ is an open continuous homomorphism of a $P\mathcal{M}$-factorizable topological group onto a topological group $H$, then $H$ is $P\mathcal{M}$-factorizable.
\end{theorem}

\begin{proof}
Indeed, it was proved in \cite[Corollary 3.8]{ZH2020} that a quotient group of a $\mathcal{M}$-factorizable topological group is also $\mathcal{M}$-factorizable. Moreover, by \cite[Corollary 4.3.24]{AA}, if $f:G\rightarrow H$ is an open continuous homomorphism of a feathered topological group onto a topological group $H$, then $H$ is also feathered. Therefore, if $G$ is $P\mathcal{M}$-factorizable, that is, feathered and $\mathcal{M}$-factorizable by Theorem \ref{PM=M}, then the topological group $H$ is also $P\mathcal{M}$-factorizable as an open continuous homomorphic image.
\end{proof}

From \cite[Proposition 3.4.2]{AA}, if a topological group $H$ is a continuous homomorphic image of an $\omega$-narrow topological group $G$, then $H$ is also $\omega$-narrow. The following corollary is follows from Theorem \ref{omega-PM}.

\begin{corollary}\cite[Theorem 2.9]{Peng2020}\label{open+continuous}
If $f:G\rightarrow H$ is an open continuous homomorphism of a $P\mathbb{R}$-factorizable topological group onto a topological group $H$, then $H$ is $P\mathbb{R}$-factorizable.
\end{corollary}

\section{Products of $P\mathcal{M}$-factorizable topological groups}

In this section, we investigate some properties about products of $P\mathcal{M}$-factorizable topological groups. In particular, some interesting properties of $\mathcal{M}$-factorizable topological groups in \cite{HW2021, ZH2020} are strengthened to $P\mathcal{M}$-factorizable topological groups. For axample, the product $G=\prod_{n\in \mathbb{N}}G_{n}$ of countably many $P\mathcal{M}$-factorizable topological groups is $P\mathcal{M}$-factorizable if and only if every factor $G_{n}$ is metrizable or every $G_{n}$ is $P\mathbb{R}$-factorizable, the product of a $P\mathcal{M}$-factorizable topological group with a compact metrizable topological group is $P\mathcal{M}$-factorizable.

First, according to the result that a locally compact group $G$ is $\mathcal{M}$-factorizable if and only if $G$ is metrizable or $G$ is $\sigma$-compact, H. Zhang et al. gave an example to show that a product of two $\mathcal{M}$-factorizable topological groups may fail to be $\mathcal{M}$-factorizable. By further observation about the example, we find that a product of two $P\mathcal{M}$-factorizable topological groups may not be $\mathcal{M}$-factorizable, so naturally not be $P\mathcal{M}$-factorizable.

\begin{example}
Assume that $G$ is a metrizable locally compact group which is not $\sigma$-compact and $H$ is a compact and non-metrizable group. Obviously, both $G$ and $H$ are $\mathcal{M}$-factorizable. Moreover, each locally compact topological group is feathered, then $G$ and $H$ both are $P\mathcal{M}$-factorizable. However, the product group $G\times H$ is neither metrizable nor $\sigma$-compact. Therefore, the product group $G\times H$ is not $P\mathcal{M}$-factorizable by Proposition \ref{3mt10}. ( Indeed, since $G\times H$ is feathered but not $\mathcal{M}$-factorizable, it can also be yielded that it is not $P\mathcal{M}$-factorizable by Theorem \ref{PM=M}.)
\end{example}

\begin{theorem}\label{product}
The product $G=\prod_{n\in \mathbb{N}}G_{n}$ of countably many $P\mathcal{M}$-factorizable topological groups is $P\mathcal{M}$-factorizable if and only if every factor $G_{n}$ is metrizable or every $G_{n}$ is $P\mathbb{R}$-factorizable.
\end{theorem}

\begin{proof}
It follows from \cite[Proposition 4.3.13]{AA} that the product space $G$ is feathered. If $G$ is $P\mathcal{M}$-factorizable, by Theorem \ref{metrizable-PR}, $G$ is either metrizable or $P\mathbb{R}$-factorizable. If $G$ is metrizable, each $G_{n}$ is also metrizable. If $G$ is $P\mathbb{R}$-factorizable, so is every $G_{n}$ by Corollary \ref{open+continuous}.

On the contrary, if every $G_{n}$ is metrizable, it is clear that $G$ is also metrizable. If every $G_{n}$ is a $P\mathbb{R}$-factorizable topological group, it is easy to see that $G$ is $\mathbb{R}$-factorizable. Moreover, $G$ is feathered, we conclude that $G$ is $P\mathbb{R}$-factorizable by Corollary \ref{PR=R}.
\end{proof}

The product of countably many $P\mathbb{R}$-factorizable topological groups is also $P\mathbb{R}$-factorizable, see \cite[Proposition 2.7]{Peng2020}, then it is clear to achieve the following by Theorems \ref{metrizable-PR} and \ref{product}.

\begin{proposition}
If $G$ is a $P\mathcal{M}$-factorizable topological group, then so is $G^{\omega}$.
\end{proposition}

\begin{remark}
Let $G$ be a compact group with $w(G)>\omega$ and $D$ an uncountable discrete group. Since $G$ and $D$ both are feathered, it is clear that $G\times D$ is feathered. However, $G$ is not metrizable and $D$ is not $P\mathbb{R}$-factorizable, then $G\times D$ is not $P\mathcal{M}$-factorizable, hence is also not $\mathcal{M}$-factorizable by Theorem \ref{PM=M}.
\end{remark}

By Lemma 3.1 and Theorem 4.7 of \cite{HW2021}, an $\mathcal{M}$-factorizable topological group which contains a non-metrizable pseudocompact subspace is $\omega$-narrow.

\begin{theorem}
Let $G$ and $H$ be topological groups, where $G$ contains a non-metrizable pseudocompact subspace. If $G\times H$ is $P\mathcal{M}$-factorizable, then $G\times H$ is $P\mathbb{R}$-factorizable.
\end{theorem}

\begin{proof}
Since every $P\mathcal{M}$-factorizable topological group is $\mathcal{M}$-factorizable, it is clear that $H$ is $\omega$-narrow by \cite[Theorem 5.7]{HW2021}. Since projection is an open continuous homomorphism, $G$ and $H$ both are $P\mathcal{M}$-factorizable by Theorem \ref{PMopen+continuous}. Then $G$ is $\omega$-narrow since $G$ contains a non-metrizable pseudocompact subspace. We have that the $P\mathcal{M}$-factorizable topological group $G\times H$ is also $\omega$-narrow, hence is $P\mathbb{R}$-factorizable by Theorem \ref{omega-PM}.
\end{proof}

\begin{theorem}
Let $G$ be a feathered group and $K$ a pseudocompact feathered group. Then $G\times K$ is $P\mathcal{M}$-factorizable if and only if either both $G$ and $K$ are metrizable or $G$ is $P\mathbb{R}$-factorizable.
\end{theorem}

\begin{proof}
Let the product group $G\times K$ be $P\mathcal{M}$-factorizable. Then the factors $G$ and $K$ are $P\mathcal{M}$-factorizable as the open continuous images by Theorem \ref{PMopen+continuous}. If $G\times K$ is not metrizable, then either $G$ or $K$ is not metrizable. If $G$ is not metrizable, it follows from Theorem \ref{metrizable-PR} that $G$ is $P\mathbb{R}$-factorizable. On the other case, if $K$ is not metrizable, $K$ is a non-metrizable pseudocompact topological group, then $G$ is $\omega$-narrow. By \cite[Theorem 2.6]{Peng2020}, a topological group is feathered and $\omega$-narrow if and only if it is $P\mathbb{R}$-factorizable, hence $G$ is a $P\mathbb{R}$-factorizable topological group.

On the contrary, if $G$ and $K$ are metrizable topological groups, it is clear that $G\times K$ is $P\mathcal{M}$-factorizable. If $G$ is $P\mathbb{R}$-factorizable, then $G\times K$ is $\mathbb{R}$-factorizable as $K$ is pseudocompact. Moreover, both $G$ and $K$ are feathered, then $G\times K$ is also feathered, which deduces that $G\times K$ is a $P\mathbb{R}$-factorizable topological group by Corollary \ref{PR=R}, hence is $P\mathcal{M}$-factorizable.
\end{proof}

\begin{theorem}\label{4dingli1}
Let $G$ and $H$ be topological groups, where the Ra\v\i kov completion $\varrho G$ of $G$ contains a non-metrizable compact subspace. If $G\times H$ is $P\mathcal{M}$-factorizable, then $H$ is $P\mathbb{R}$-factorizable.
\end{theorem}

\begin{proof}
Since every $P\mathcal{M}$-factorizable topological group is $\mathcal{M}$-factorizable, it follows from \cite[Theorem 3.11]{HW2021-2} that $H$ is pseudo-$\omega_{1}$-compact, hence $H$ is $\omega$-narrow by \cite[Proposition 3.4.31]{AA}. Since the product $G\times H$ is $P\mathcal{M}$-factorizable, so is $H$ by Theorem \ref{PMopen+continuous}. Therefore, we have that $H$ is a $P\mathbb{R}$-factorizable topological group by Theorem \ref{omega-PM}.
\end{proof}

\begin{proposition}
If the product $G\times H$ of topological groups is $P\mathcal{M}$-factorizable and the group $G$ is precompact, then either $G$ is second countable or $H$ is $P\mathbb{R}$-factorizable.
\end{proposition}

\begin{proof}
The first part that $G$ is second countable follows just from \cite[Proposition 3.12]{HW2021-2} and the second part is deduced by Theorem \ref{4dingli1}.
\end{proof}

\begin{theorem}
Let $G$ be a feathered group and $H$ a precompact feathered group. Then $G\times H$ is $P\mathcal{M}$-factorizable if and only if either both $G$ and $H$ are metrizable or $G$ is Lindel\"{o}f $\Sigma$-group.
\end{theorem}

\begin{proof}
The necessity is claimed in \cite[Theorem 3.13]{HW2021-2}, where $H$ just need to be precompact.

It suffices to prove the sufficiency. On the first case, if both $G$ and $H$ are metrizable, it is trivial that $G\times H$ is $P\mathcal{M}$-factorizable. On the other case, let $G$ be a feathered Lindel\"{o}f $\Sigma$-group and $H$ a precompact feathered group. Then the Ra\v\i kov completion $\varrho H$ of $H$ is compact. $G\times H$ is $\mathbb{R}$-factorizable as a subgroup of the Lindel\"{o}f $\Sigma$-group $G\times \varrho H$. Moreover, since both $G$ and $H$ are feathered, $G\times H$ is also feathered, hence is $P\mathbb{R}$-factorizable by Corollary \ref{PR=R}. We obtain that $G\times H$ is a $P\mathcal{M}$-factorizable topological group.
\end{proof}

W. He et al. proved that the product of an $\mathcal{M}$-factorizable topological group with a locally compact separable metrizable topological group is $\mathcal{M}$-factorizable, see \cite[Theorem 3.14]{HW2021-2}.

\begin{theorem}\label{local+separ+metri}
Let $G$ be a $P\mathcal{M}$-factorizable topological group and $H$ a locally compact separable metrizable topological group. Then $G\times H$ is $P\mathcal{M}$-factorizable.
\end{theorem}

\begin{proof}
Since the product of an $\mathcal{M}$-factorizable topological group with a locally compact separable metrizable topological group is $\mathcal{M}$-factorizable, $G\times H$ is $\mathcal{M}$-factorizable. Moreover, both $G$ and $H$ are feathered as every locally compact group is feathered, so $G\times H$ is also feathered. Therefore, $G\times H$ is $P\mathcal{M}$-factorizable by Theorem \ref{PM=M}.
\end{proof}

\begin{remark}
Indeed, by revising the proof to \cite[Theorem 3.14]{HW2021-2}, we can give a direct proof of Theorem \ref{local+separ+metri}
\end{remark}

\begin{proof}
Since $H$ is a locally compact separable metrizable topological group, $H$ is $\sigma$-compact. Then there is an increasing sequence $\{H_{n}:n\in \mathbb{N}\}$ of compact subsets of $H$ such that $H=\bigcup_{n\in \mathbb{N}}H_{n}$ and $H_{n}$ is contained in the interior of $H_{n+1}$ for each $n\in \mathbb{N}$. Let $f$ be a continuous real-valued function on $G\times H$. Denote by $C(H_{n})$ the space of continuous real-valued functions on $H_{n}$  with sup-norm, for each $n\in \mathbb{N}$. Then define a mapping $\Psi_{n}:G\rightarrow C(H_{n})$ by $\Psi_{n}(x)(y)=f(x,y)$ for all $x\in G$ and $y\in H_{n}$. Since $H_{n}$ is compact and second countable, $\Psi_{n}$ is continuous and $C(H_{n})$ is second countable. Put $\Psi$ the diagonal product of $\{\Psi_{n}:n\in \mathbb{N}\}$. Since $\prod_{n\in \mathbb{N}}C(H_{n})$ is second countable, it is clear that $\Psi (G)$ is also second countable.

By the hypothesis, $G$ is $P\mathcal{M}$-factorizable, we can find a perfect homomorphism $\pi$ of $G$ onto a metrizable group $K$ and a  continuous mapping $\psi$ of $K$ to $\Psi (G)$ such that $\Psi =\psi \circ \pi$. Take $y\in H$ and choose $n\in \mathbb{N}$ with $y\in H_{n}$. Let $x,x'\in G$. If $\pi (x)=\pi (x')$, then $\Psi (x)=\Psi (x')$. Then $\Psi_{n}(x)=\Psi_{n}(x')$, that is, $f(x,y)=f(x',y)$. Therefore, we can define a mapping $h:K\times H\rightarrow \mathbb{R}$ such that $f=h\circ (\pi \times id_{H})$. It is not difficult to verify that $h$ is continuous. Since both $K$ and $H$ are metrizable topological groups, $G\times H$ is also metrizable. Moreover, $\pi$ is a perfect mapping, so is the mapping $\pi \times id_{H}$. Thus, we conclude that the product $G\times H$ is $P\mathcal{M}$-factorizable.
\end{proof}

\begin{theorem}
Let $G$ be a $P\mathbb{R}$-factorizable topological group and $H$ a locally compact separable metrizable topological group. Then $G\times H$ is $P\mathbb{R}$-factorizable.
\end{theorem}

\begin{proof}
First, $G\times H$ is $P\mathcal{M}$-factorizable by Theorem \ref{local+separ+metri}. Moreover, since both $G$ and $H$ are $\omega$-narrow, then $G\times H$ is $\omega$-narrow and it is achieved that it is $P\mathbb{R}$-factorizable by Theorem \ref{omega-PM}.
\end{proof}

\begin{corollary}\label{PM+com-metri}
Let $G$ be a $P\mathcal{M}$-factorizable topological group and $H$ a compact metrizable topological group. Then $G\times H$ is $P\mathcal{M}$-factorizable.
\end{corollary}

\begin{corollary}\label{PR*compact+metrizable}
Let $G$ be a $P\mathbb{R}$-factorizable topological group and $H$ a compact metrizable topological group. Then $G\times H$ is $P\mathbb{R}$-factorizable.
\end{corollary}

A topological space $X$ is called {\it pseudo-$\aleph_{1}$-compact} if every discrete family of open subsets of $X$ is countable. As we all know, every separable space is pseudo-$\aleph_{1}$-compact and every pseudo-$\aleph_{1}$-compact topological group is $\omega$-narrow. Therefore, the following resut follows from Theorem \ref{omega-PM}.

\begin{corollary}
A pseudo-$\aleph_{1}$-compact $P\mathcal{M}$-factorizable topological group is $P\mathbb{R}$-factorizable.
\end{corollary}

\begin{corollary}\label{pseudo+PM}
A product of a pseudo-$\aleph_{1}$-compact $P\mathcal{M}$-factorizable topological group and a compact group is $P\mathbb{R}$-factorizable.
\end{corollary}

\begin{proof}
Indeed, let $G$ be a pseudo-$\aleph_{1}$-compact $P\mathcal{M}$-factorizable topological group and $H$ a compact group. It was proved in \cite[Corollary 5.2]{ZH2020} that a product of a pseudo-$\aleph_{1}$-compact $\mathcal{M}$-factorizable topological group and a compact group is $\mathbb{R}$-factorizable. Then the product $G\times H$ is $\mathbb{R}$-factorizable. Moreover, $G$ is feathered by Proposition \ref{PM-feathered} and $H$ is also feathered, so $G\times H$ is a feathered group, which deduces that $G\times H$ is $P\mathbb{R}$-factorizable by Corollary \ref{PR=R}.
\end{proof}

It was proved in \cite[Theorem 5.4]{ZH2020} that if $G\times K$ is $\mathcal{M}$-factorizable, where $G$ is an $\mathcal{M}$-factorizable group and $K$ is compact group, then $G$ is pseudo-$\aleph_{1}$-compact or $K$ is metrizable. Then we have the following by Corollaries \ref{PM+com-metri} and \ref{pseudo+PM}.

\begin{theorem}\label{PM*compact}
Let $G$ be a $P\mathcal{M}$-factorizable topological group and $K$ a compact group. Then $G\times K$ is $P\mathcal{M}$-factorizable if and only if one of the following conditions holds:
\begin{enumerate}
\smallskip
\item $K$ is metrizable;

\smallskip
\item $G$ is pseudo-$\aleph_{1}$-compact.
\end{enumerate}
\end{theorem}

Recall that a mapping $f:X\rightarrow Y$ is {\it $d$-open} if for every open set $U$ in $X$, the image $f(U)$ is contained in the interior of the closure of $f(U)$. The following results was proved in \cite{HW2021}, see Proposition 6.3 and Theorem 6.5.

\begin{proposition}
An image of a feathered topological group under a continuous $d$-open homomorphism is feathered.
\end{proposition}

\begin{proposition}
Let $p$ be a continuous $d$-open homomorphism from a topological group $G$ onto a topological group $H$. If $G$ is $\mathcal{M}$-factorizable, so is $H$.
\end{proposition}

Since it is proved in Theorem \ref{PM=M} that a topological group $G$ is $P\mathcal{M}$-factorizable if and only if $G$ is feathered $\mathcal{M}$-factorizable, the following is deduced by two propositions above.

\begin{corollary}
If a topological group $H$ is a continuous $d$-open homomorphic image of a $P\mathcal{M}$-factorizable topological group $G$, then $H$ is also $P\mathcal{M}$-factorizable.
\end{corollary}

\section{On $Pm$-factorizable topological groups}

In this section, according to the concept of $m$-factorizable topological groups, that is, a topological group $G$ is called {\it $m$-factorizable} if for every continuous mapping $f:G\rightarrow M$ to a metrizable space $M$, there exists a continuous homomorphism $\pi :G\rightarrow K$ onto a second-countable group $K$ such that $f=g\circ \pi$, for some continuous homomorphism $g$ from $K$ onto $M$, we introduce $Pm$-factorizable topological groups  as the following by strengthening the continuous homomorphism $\pi$ to a perfect homomorphism.

\begin{definition}\cite{AA}
A topological group $G$ is called {\it $m$-factorizable} if for every continuous mapping $f:G\rightarrow M$ to a metrizable space $M$, there exists a continuous homomorphism $\pi :G\rightarrow K$ onto a second-countable group $K$ such that $f=g\circ \pi$, for some continuous homomorphism $g$ from $K$ onto $M$.
\end{definition}

In particular, if the continuous homomorphism $\pi$ is perfect, we call the topological group $G$ {\it $Pm$-factorizable}. Clearly, every $Pm$-factorizable topological group is $m$-factorizable and every $Pm$-factorizable topological group is $P\mathbb{R}$-factorizable. Then we show that a topological group $G$ is $Pm$-factorizable if and only if $G$ is $P\mathbb{R}$-factorizable and pseudo-$\aleph_{1}$-compact and the product of a $Pm$-factorizable topological group with an arbitrary compact group is $Pm$-factorizable.

First, by \cite[Proposition 8.5.1]{AA}, we just need to revise the continuous homomorphism $\pi$ to a perfect homomorphism $\pi$, we can obtain the following.

\begin{proposition}
A topological group $G$ is $Pm$-factorizable if and only if for every continuous pseudometric $d$ on $G$, there exist a perfect homomorphism $\pi:G\rightarrow K$ onto a second-countable topological group $K$ and a continuous pseudometric $\varrho $ on $K$ such that $d(x,y)=\varrho (\pi (x),\pi (y))$, for all $x,y\in G$.
\end{proposition}

\begin{lemma}\cite[Lemma 3.2]{Peng2020}\label{8.1.2-perfect}
Suppose that $f:G\rightarrow X$ is a continuous mapping of a $P\mathbb{R}$-factorizable topological group $G$ to a Tychonoff space $X$ with $w(X)\leq \tau$. Then one can find a perfect homomorphism $\pi:G\rightarrow K$ onto a topological group $K$ with $w(K)\leq \tau$ such that $f=h\circ g$ for some continuous mapping $h:g(G)\rightarrow K$.
\end{lemma}

\begin{theorem}\label{PR+pseudo}
A topological group $G$ is $Pm$-factorizable if and only if $G$ is $P\mathbb{R}$-factorizable and pseudo-$\aleph_{1}$-compact.
\end{theorem}

\begin{proof}
Since every $Pm$-factorizable topological group is $m$-factorizable and every $Pm$-factorizable topological group is $P\mathbb{R}$-factorizable, the necessity is clear since all $m$-factorizable topological groups are pseudo-$\aleph_{1}$-compact. It suffices to claim the sufficiency.

Assume that a topological group $G$ is $P\mathbb{R}$-factorizable and pseudo-$\aleph_{1}$-compact and $f:G\rightarrow M$ is a continuous mapping of $G$ onto a metrizable space $M$. Since continuous mapping preserve pseudo-$\aleph_{1}$-compactness, $M$ is also pseudo-$\aleph_{1}$-compact, then $w(M)\leq \omega$. Since $G$ is $P\mathbb{R}$-factorizable, it follows from Lemma \ref{8.1.2-perfect} that we can find a perfect homomorphism $\pi :G\rightarrow K$ onto a second-countable topological group $K$ such that $f=g\circ \pi$ for some continuous mapping $g:K\rightarrow M$ onto the metrizable space $M$. Therefore, we have that $G$ is a $Pm$-factorizable topological group.
\end{proof}

Since a topological group $G$ is $m$-factorizable if and only if $G$ is $\mathbb{R}$-factorizable and pseudo-$\aleph_{1}$-compact by \cite[Theorem 8.5.2]{AA}, it is clear to deduce the following by Corollary \ref{PR=R}.

\begin{corollary}\label{Pm=m}
A topological group $G$ is $Pm$-factorizable if and only if $G$ is feathered and $m$-factorizable.
\end{corollary}

The following corollary is from Theorem \ref{PR+pseudo} and Corollary \ref{open+continuous}

\begin{corollary}\label{Pmopen+continuous}
Let $\pi:G\rightarrow H$ be an open continuous homomorphism of a topological group $G$ onto $H$. If $G$ is $Pm$-factorizable, so is $H$.
\end{corollary}

By \cite[Theorem 8.5.8]{AA}, if an $\mathbb{R}$-factorizable group $G$ satisfies $w(G)\leq \tau \geq \aleph_{0}$, then $|C(G)|\leq \tau^{\omega}$, where $C(X)$ denotes the set of continuous real-valued functions on a space $X$. Indeed, if $G$ is an $\mathbb{R}$-factorizable group with $w(G)^{\omega}<2^{\aleph_{1}}$, then $G$ is pseudo-$\aleph_{1}$-compact. Then the following result follows from Theorem \ref{PR+pseudo}.

\begin{corollary}
Every $P\mathbb{R}$-factorizable group $G$ with $w(G)^{\omega}<2^{\aleph_{1}}$ is $Pm$-factorizable.
\end{corollary}

\begin{corollary}
Let $G$ be a $P\mathbb{R}$-factorizable group with $w(G)^{\omega}<2^{\aleph_{1}}$. Then the product $G\times K$ is $Pm$-factorizable, for every compact group $K$.
\end{corollary}

\begin{lemma}\label{lemma 5.6}
Let $G$ be a topological group with the property that for every continuous function $f:G\rightarrow \mathbb{R}$, there exists a perfect homomorphism $\pi :G\rightarrow H$ onto a $P\mathbb{R}$-factorizable group $H$ such that $f=g\circ \pi$ for some continuous function $g:H\rightarrow \mathbb{R}$. Then the group $G$ is $P\mathbb{R}$-factorizable.
\end{lemma}

\begin{proof}
Let $f:G\rightarrow \mathbb{R}$ be a continuous function. By the hypothesis, there exists a perfect homomorphism $\pi :G\rightarrow H$ onto a $P\mathbb{R}$-factorizable group $H$ and a continuous function $g:H\rightarrow \mathbb{R}$ such that $f=g\circ \pi$. Since $H$ is $P\mathbb{R}$-factorizable, we can find a perfect homomorphism $p:H\rightarrow K$ onto a second-countable topological group $K$ and a continuous real-valued function $h$ on $K$ such that $g=h\circ p$.

\begin{equation*}
\centerline{
\xymatrix{
 G \ar[r]^{f } \ar[d]_{\pi} & \mathbb{R} \\
 H \ar[ur]^{g} \ar[r]^{p}&   K \ar[u]^{h}   & &  }}
\end{equation*}
Let $\varphi =p\circ \pi$. Since both $p$ and $\pi$ are perfect, $\varphi$ is also perfect. Therefore, there is a perfect homomorphism $\varphi$ of $G$ onto $K$ satisfying $f=h\circ \varphi$. We conclude that $G$ is a $P\mathbb{R}$-factorizable topological group.
\end{proof}

\begin{theorem}\label{Pm*compact}
Let $G$ be a $Pm$-factorizable topological group and $K$ an arbitrary compact group, then the topological group $G\times K$ is $Pm$-factorizable.
\end{theorem}

\begin{proof}
By Corollary \ref{Pm=m}, $G$ is feathered and $m$-factorizable. Since the product group of every $m$-factorizable group with an arbitrary compact group is also $m$-factorizable, it is clear that $G\times K$ is $m$-factorizable. Moreover, every $Pm$-factorizable topological group is feathered and every locally compact group is also feathered, so $G\times K$ is feathered, which deduces that $G\times K$ is $Pm$-factorizable By Corollary \ref{Pm=m}.
\end{proof}

\begin{remark}
The following is a direct proof of Theorem \ref{Pm*compact}.
\end{remark}

\begin{proof}
Let $f:G\times K\rightarrow \mathbb{R}$ be a continuous function. $C(K)$ denotes the space of all continuous real-valued functions on $K$ with the sup-norm topology and consider the mapping $\Psi :G\rightarrow C(K)$ defined by $\Psi (x)(y)=f(x,y)$, for all $x\in G$ and $y\in K$. Since $K$ is compact, $\Psi$ is continuous. Since $G$ is a $Pm$-factorizable topological group, it follows from Theorem \ref{PR+pseudo} that $G$ is $P\mathbb{R}$-factorizable and pseudo-$\aleph_{1}$-compact, so the subspace $\Psi (G)$ of the metric space $C(K)$ is pseudo-$\aleph_{1}$-compact and hence is second-countable. Since $G$ is $P\mathbb{R}$-factorizable,  we can find a perfect homomorphism $\pi :G\rightarrow H$ onto a second-countable topological group $H$ and a continuous mapping $\psi :H\rightarrow C(K)$ such that $\Psi =\psi\circ\pi$ by Lemma \ref{8.1.2-perfect}.

Then if $x_{1},x_{2}\in G$ and $\pi (x_{1})=\pi (x_{2})$, then $f(x_{1},y)=f(x_{2},y)$ for each $y\in K$. Suppose on the contrary, if $f(x_{1},y)\not=f(x_{2},y)$ for some $x_{1},x_{2}\in G$ and $y\in K$, then $\Psi (x_{1})(y)=\Psi (x_{2})(y)$. It is easy to see that $\pi (x_{1})\not=\pi (x_{2})$. Therefore, we can define a mapping $h:H\times K\rightarrow \mathbb{R}$ such that $h\circ (\pi\times id_{K})=f$, where $id_{K}$ is the identity mapping of $K$ onto itself.
\begin{equation*}
\centerline{
\xymatrix{
  G\times K \ar[d]_{\pi \times id_{K}} \ar[r]^{f} & \mathbb{R}       \\
  H\times K \ar[ur]_{h}                     }}
\end{equation*}

Moreover, $H$ is continuous. For an arbitrary point $(s,y)\in H\times K$ and a number $\varepsilon >0$. Let $x^{*}\in G$ be such that $\pi (x^{*})=s$. Since $\psi$ is continuous, there exists an open neighborhood $U$ of $s$ in $H$ such that $||\Psi (t)-\Psi (s)||<\varepsilon /2$ for each $t\in U$. We can find a neighborhood $V$ of $y$ in $K$ such that $|f(x^{*},z)-f(x^{*},y)|<\varepsilon /2$ for each $z\in V$. Let $(t,z)\in U\times V$  and $x\in G$ such that $\pi (x)=t$. Then $$|h(t,z)-h(s,y)|\leq |f(x,z)-f(x^{*},z)|+|f(x^{*},z)-f(x^{*},y)|<\varepsilon /2+\varepsilon /2=\varepsilon.$$ Therefore, $h$ is continuous.

Since $H$ is a second-countable topological group and $K$ is a compact group, we know that $H\times K$ is Lindel\"{o}f. Moreover, it is well-known that every metrizable topological group is feathered and every locally compact topological group is feathered, it follows that $H\times K$ is feathered. As every Lindel\"{o}f topological group is $\mathbb{R}$-factorizable, we obtain that $H\times K$ is $P\mathbb{R}$-factorizable by Corollary \ref{PR=R}. Since $f=h\circ \pi\times id_{K}$, it follows from Lemma \ref{lemma 5.6} that $G\times K$ is $P\mathbb{R}$-factorizable. Finally, as $G$ is pseudo-$\aleph_{1}$-compact and $K$ is compact, the product $G\times K$ is pseudo-$\aleph_{1}$-compact and we conclude that $G\times K$ is $Pm$-factorizable by Theorem \ref{PR+pseudo}.
\end{proof}

In \cite[Theorem 8.5.11]{AA}, if the product $G\times K$ of an $\mathbb{R}$-factorizable group $G$ and a compact group $K$ is $\mathbb{R}$-factorizable, then either $G$ is pseudo-$\aleph_{1}$-compact or $K$ is metrizable. Therefore, by Theorem \ref{Pm*compact} and Corollary \ref{PR*compact+metrizable}, we have the following result which is similar to Theorem \ref{PM*compact}.

\begin{theorem}\label{PR*compact}
Let $G$ be a $P\mathbb{R}$-factorizable topological group and $K$ a compact group. Then $G\times K$ is $P\mathbb{R}$-factorizable if and only if one of the following conditions holds:
\begin{enumerate}
\smallskip
\item $K$ is metrizable;

\smallskip
\item $G$ is pseudo-$\aleph_{1}$-compact.
\end{enumerate}
\end{theorem}

\begin{proposition}
A $C$-embedded closed subgroup of a $P\mathbb{R}$-factorizable topological group is $P\mathbb{R}$-factorizable.
\end{proposition}

\begin{proof}
Let $H$ be a $C$-embedded closed subgroup of a $P\mathbb{R}$-factorizable topological group $G$. Then every continuous function $f:H\rightarrow \mathbb{R}$ admits an extension to a continuous function $g:G\rightarrow \mathbb{R}$. Since $G$ is $P\mathbb{R}$-factorizable, there exists a perfect homomorphism $\pi :G\rightarrow K$ onto a second-countable group $K$ such that $g=h\circ \pi$ for some continuous function $h$ on $K$. Since $H$ is a closed subgroup of $G$, the homomorphism $\pi\!\!\upharpoonright_{H}$ of $H$ onto the subgroup $\pi (H)$ of $K$ is also perfect and factorizes $f$.
\end{proof}

The following follows from the fact that every retract of a space $X$ is $C$-embedded in $X$.

\begin{corollary}\label{retract}
Let $G$ be a $P\mathbb{R}$-factorizable topological group and $H$ a closed subgroup of $G$. If $H$ is a retract of $G$, then $H$ is also $P\mathbb{R}$-factorizable.
\end{corollary}

\begin{corollary}
Let $G$ be a topological group. If $G\times \mathbb{Z}(2)^{\omega_{1}}$ is $P\mathbb{R}$-factorizable, then $G$ is $Pm$-factorizable.
\end{corollary}

\begin{proof}
Let $K=\mathbb{Z}(2)^{\omega_{1}}$, where $\mathbb{Z}(2)=\{0,1\}$ is the discrete group. Indeed, it was proved in \cite[Theorem 8.5.5]{AA} that if $G\times \mathbb{Z}(2)^{\omega_{1}}$ is $\mathbb{R}$-factorizable, $G$ is pseudo-$\aleph_{1}$-compact. Then it suffices to show that $G$ is $P\mathbb{R}$-factorizable by Theorem \ref{PR+pseudo}.

Let $e_{K}$ be the identity of the group $K$. It is clear that $G\cong G\times \{e_{K}\}$ is closed in $G\times K$ and is a retract of $G\times K$. By Corollary \ref{retract}, $G$ is $P\mathbb{R}$-factorizable, which deduces that $G$ is a $Pm$-factorizable topological group.
\end{proof}

\section{On $S\mathcal{M}$-factorizable topological groups}

In this section, we pose the concepts of $S\mathcal{M}$-factorizable topological groups and $PS\mathcal{M}$-factorizable topological groups. We show that a topological group $G$ is $PS\mathcal{M}$-factorizable if and only if $G$ is feathered $S\mathcal{M}$-factorizable, and the properties of $S\mathcal{M}$-factorizabilities and $PS\mathcal{M}$-factorizabilities are preserved by open continuous homomorphisms on topological groups. Indeed, by the definition of $m$-factorizable topological groups, it is natural to extent $\mathbb{R}$ to a metrizable space, so we pose the following concept.

\begin{definition}
A topological group $G$ is called {\it strongly $\mathcal{M}$-factorizable} ($S\mathcal{M}$-factorizable for short) if for every continuous mapping $f:G\rightarrow M$ to a metrizable space $M$, there exists a continuous homomorphism $\pi :G\rightarrow H$ onto a metrizable group $H$ and a continuous mapping $g:H\rightarrow M$ such that $f=g\circ \pi$. In particular, if the continuous homomorphism $\pi :G\rightarrow H$ is perfect, we call $G$ {\it $PS\mathcal{M}$-factorizable.}
\end{definition}

\smallskip
\begin{equation*}
\centerline{
\xymatrix{
  G \ar[d]_{\pi } \ar[r]^{f} & M       \\
  H \ar[ur]_{g}                     }}
\end{equation*}

Clearly, every $S\mathcal{M}$-factorizable topological group is $\mathcal{M}$-factorizable and every $m$-factorizable topological group is $S\mathcal{M}$-factorizable.

\begin{theorem}\label{omega+SM}
A topological group $G$ is $m$-factorizable if and only if $G$ is $S\mathcal{M}$-factorizable and $\omega$-narrow.
\end{theorem}

\begin{proof}
Since every $m$-factorizable topological group is $\mathbb{R}$-factorizable, hence is $\omega$-narrow, the necessity is trivial.

Let's prove the sufficiency. Suppose that $G$ is an $S\mathcal{M}$-factorizable and $\omega$-narrow topological group, $f:G\rightarrow M$ is a continuous mapping to a metrizable space $M$. Then there exists a continuous homomorphism $\varphi :G\rightarrow K$ onto a metrizable topological group $K$ such that $f=g\circ \varphi$, where $g:K\rightarrow M$ is a continuous mapping. Since $G$ is $\omega$-narrow, so is $K$. Since every first-countable $\omega$-narrow topological group is second-countable, we obtain that $G$ is $m$-factorizable.
\end{proof}

Similarly, the following result also holds.

\begin{proposition}
A topological group $G$ is $Pm$-factorizable if and only if $G$ is $PS\mathcal{M}$-factorizable and $\omega$-narrow.
\end{proposition}

\begin{proposition}\label{SM===}
A topological group $G$ is $S\mathcal{M}$-factorizable if and only if for every continuous pseudometric $d$ on $G$, one can find a continuous homomorphism $\pi:G\rightarrow K$ onto a metrizable topological group $K$ and a continuous pseudometric $\varrho$ on $K$ such that $d(x,y)=\varrho (\pi (x),\pi (y))$ for all $x,y\in G$.
\end{proposition}

\begin{proof}
Let $G$ be an $S\mathcal{M}$-factorizable topological group and $d$ a continuous pseudometric on $G$. Consider the metric space $M=G/d$ with the associated metric $d^{*}$, obtained from $G$ by identifying points at zero distance with respect to $d$. Let $p:G\rightarrow G/d$ be the projection assigning to a point $x\in G$ the equivalence class $\overline{x}$ consisting of all $z\in G$ with $d(x,z)=0$. Then $d(x,y)=d^{*}(p(x),p(y))$, for all $x,y\in G$. By the hypothesis, $G$ is $S\mathcal{M}$-factorizable, there exists a continuous homomorphism $\pi :G\rightarrow K$ onto a metrizable topological group $K$ and a continuous mapping $h:K\rightarrow M$ such that $p=h\circ \pi$. Then define a continuous pseudometric $\varrho$ on $K$ by $\varrho (s,t)=d^{*}(h(s),h(t))$ for all $s,t\in K$. Then for all $x,y\in G$, $$d(x,y)=d^{*}(p(x),p(y))=d^{*}(h(\pi (x)),h(\pi (y)))=\varrho (\pi (x),\pi (y)).$$

On the other hand, let $f:G\rightarrow M$ be a continuous mapping to a metric space $M$ with a metric $\kappa$. Define a continuous pseudometric $d$ on $G$ by $d(x,y)=\kappa (f(x),f(y))$ for all $x,y\in G$. By the hypothesis, there exists a continuous homomorphism $\pi :G\rightarrow K$ onto a metrizable topological group $K$ and a continuous pseudometric $\varrho$ on $K$ such that $d(x,y)=\varrho (\pi (x),\pi (y))$, for all $x,y\in G$. Therefore, $\varrho (\pi (x),\pi (y))=\kappa (f(x),f(y))$, for all $x,y\in G$. Then $\pi (x)=\pi (y)$ always implies $f(x)=f(y)$. We can find a continuous mapping $h:K\rightarrow M$ such that $f=h\circ \pi$. Thus, we conclude that $G$ is $S\mathcal{M}$-factorizable.
\end{proof}

By Proposition \ref{SM===}, it is not difficult to see that if we change the continuous homomorphism $\pi$ onto a perfect homomorphism, it also holds for $PS\mathcal{M}$-factorizable topological groups. Then we show that every feathered $S\mathcal{M}$-factorizable is $PS\mathcal{M}$-factorizable.

\begin{theorem}\label{PSM=feathered+SM}
A topological group $G$ is $PS\mathcal{M}$-factorizable if and only if $G$ is feathered $S\mathcal{M}$-factorizable.
\end{theorem}

\begin{proof}
The necessity is trivial, it suffices to claim the sufficiency.

Let $G$ be a feathered $S\mathcal{M}$-factorizable topological group and $f:G\rightarrow H$ a continuous mapping to a metrizable space $E$. Then we can find a continuous homomorphism $\pi :G\rightarrow H$ onto a metrizable topological group $H$ and a continuous mapping $g:H\rightarrow E$ such that $f=g\circ \pi$. Since every $\mathcal{M}$-factorizable group is $\omega$-balanced by \cite[Theorem 3.1]{ZH2020} and $G$ is feathered, there exists a perfect homomorphism $p:G\rightarrow K$ onto a metrizable topological group $K$. Let $\varphi$ be the diagonal product of the homomorphisms $\pi$ and $p$ and $M=\varphi (G)\subseteq H\times K$. Since $p$ is perfect, the homomorphism $\varphi$ is also perfect by \cite[Theorem 3.7.11]{E}. By the definition of $\varphi$, we can find continuous homomorphisms $q_{H}:M\rightarrow H$ and $q_{K}:M\rightarrow K$ satisfying $\pi =q_{H}\circ\varphi$ and $p=q_{K}\circ\varphi$. Then from \cite[Proposition 3.7.5]{E}, it follows that $q_{K}$ is perfect.

Then the subgroup $M$ of $H\times K$ is metrizable. Indeed, since $H$ and $K$ are both first-countable and the property of first-countability is hereditary, it is clear that the subgroup $M$ of $H\times K$ is first-countable, hence $M$ is metrizable.

\begin{equation*}
\centerline{
\xymatrix{ & E \\
 G \ar[ur]^{f } \ar[d]_{p} \ar[dr]^{\varphi} \ar[r]^-{\pi}&  H  \ar[u]^{g} \\
 K& \ar[l]^{q_{K}}  M  \ar[u]_{q_{H}} & &  }}
\end{equation*}

We define a continuous mapping $h:M\rightarrow E$ by $h=g\circ q_{H}$. Then, for each continuous mapping $f:G\rightarrow E$ to a metrizable space $E$, we can find a perfect homomorphism $\varphi :G\rightarrow M$ onto a metrizable topological group $M$ and a continuous mapping $h:M\rightarrow E$ such that $f=h\circ\varphi$. Therefore, we conclude that $G$ is $PS\mathcal{M}$-factorizable.
\end{proof}

\begin{lemma}\label{6yl1}
Let $\{U_{n}:n\in \mathbb{N}\}$ be a family of neighborhoods of the identity $e$ in an $\omega$-balanced group $G$. Then there exists a continuous homomorphism $p:G\rightarrow H$ onto a metrizable topological group $H$ and a family $\{V_{n}:n\in \mathbb{N}\}$ of open neighborhoods of the identity $e_{H}$ in $H$ such that $p^{-1}(V_{n})\subseteq U_{n}$, for each $n\in \mathbb{N}$.
\end{lemma}

\begin{proof}
It is well-known that if $G$ is an $\omega$-balanced topological group, for each open neighborhood $U$ of the identity element $e$ in $G$, there exists a continuous homomorphism $\pi$ of $G$ onto a metrizable group $H$ such that $\pi^{-1}(V)\subseteq U$, for some open neighborhood $V$ of the identity element $e_{H}$ of $H$, see \cite[Theorem 4.3.18]{AA}. Then for every $n\in \mathbb{N}$, we can find a continuous homomorphism $\pi_{n}$ of $G$ onto a metrizable topological group $H_{n}$ and an open neighborhood $W_{n}$ of the identity element in $H_{n}$ such that $\pi_{n}^{-1}(W_{n})\subseteq U_{n}$. Let $\pi$ be the diagonal product of the family $\{\pi_{n}:n\in \mathbb{N}\}$. Then $H=\pi (G)$ is a subgroup of the product $P=\prod_{n\in \mathbb{N}}H_{n}$, so the group $H$ is also metrizable. Denote by $p_{n}$ the projection of $P$ onto the factor $H_{n}$. Then $\pi_{n}=p_{n}\circ \pi$ for each $n\in \mathbb{N}$. The open neighborhoods $V_{n}=H\cap p_{n}^{-1}(W_{n})$ of the identity in $H$ is such that $\pi^{-1}(V_{n})=\pi_{n}^{-1}(W_{n})\subseteq U_{n}$.
\end{proof}

\begin{theorem}\label{SMopen+continuous}
Let $\pi:G\rightarrow H$ be an open continuous homomorphism of a topological group $G$ onto $H$. If $G$ is $S\mathcal{M}$-factorizable, so is $H$.
\end{theorem}

\begin{proof}
Let $f:H\rightarrow M$ be a continuous mapping to a metrizable space $M$. Then $f\circ \pi :G\rightarrow M$ is a continuous mapping. Since $G$ is $S\mathcal{M}$-factorizable, there exists a continuous homomorphism $\varphi :G\rightarrow K$ onto a metrizable topological group $K$ and a continuous mapping $g:K\rightarrow M$ such that $f\circ \pi=g\circ \varphi$. Since $K$ is a first-countable topological group, we can find a countable local base $\{U_{n}:n\in \mathbb{N}\}$ at the identity $e_{K}$ of $K$ and put $V_{n}=\pi (\varphi^{-1}(U_{n}))$, for each $n\in \mathbb{N}$, $V_{n}$ is an open neighborhood of the identity $e_{H}$ in $H$, for each $n\in \mathbb{N}$.

In the metrizable space $M$, the metric which generates the original topology is denoted by $d$. For each $m\in M$ and arbitrary $\varepsilon >0$, let $O_{\varepsilon}(m)=\{t\in M:d(m,t)<\varepsilon\}$. For each $h\in H$ and $\varepsilon >0$, choose $g\in G$ with $\pi (g)=h$ and put $x=\varphi (g)$. Then $f(h)=f(\pi (g))=g(\varphi (g))=g(x)$. Since $g$ is continuous, there exists $n\in \mathbb{N}$ such that $g(xU_{n})\subseteq O_{\varepsilon}(f(h))$. Since $V_{n}=\pi (\varphi^{-1}(U_{n}))$, for each $n\in \mathbb{N}$, it is easy to see that $f(hV_{n})\subseteq O_{\varepsilon}(f(h))$ for some $n\in \mathbb{N}$.

\begin{equation*}
\centerline{
\xymatrix{ G  \ar[d]_{\varphi} \ar[r]^{\pi}&  H \ar[d]_{f}  \ar[r]^{p} &L \ar[dl]_{h} \\
 K \ar[r]^{g} & M   &  }}
\end{equation*}

Since $\pi:G\rightarrow H$ is an open continuous homomorphism of an $\omega$-balanced topological group $G$ onto topological group $H$, $H$ is also $\omega$-balanced. Then there exists a continuous homomorphism $p:H\rightarrow L$ onto a metrizable topological group $L$ with a local base $\{W_{n}:n\in \mathbb{N}\}$ at the identity such that $p^{-1}(W_{n})\subseteq V_{n}$ for each $n\in \mathbb{N}$ by Lemma \ref{6yl1}. Put $N$ the kernel of $p$. Then $N\subseteq \bigcap_{n\in \mathbb{N}}V_{n}$. Therefore, $f$ is constant on every coset $hN$ in $H$. Then we can define a mapping $h:L\rightarrow M$ such that $h\circ p=f$. It is not difficult to see that for every $y\in L$ and arbitrary $\varepsilon >0$, there exists $n\in \mathbb{N}$ such that $h(yW_{n})\subseteq O_{\varepsilon}(h(y))$, which means that $h$ is continuous. We conclude that $H$ is also $S\mathcal{M}$-factorizable.

\end{proof}

If $f:G\rightarrow H$ is an open continuous homomorphism of a feathered group $G$ onto a group $H$, $H$ is also feathered, see \cite[Corollary 4.3.24]{AA}, the following corollary is clear by Theorems \ref{PSM=feathered+SM} and \ref{SMopen+continuous}.

\begin{corollary}\label{PSMopen+continuous}
Suppose that $\pi:G\rightarrow H$ is an open continuous homomorphism of a $PS\mathcal{M}$-factorizable topological group $G$ onto $H$, then $H$ is also $PS\mathcal{M}$-factorizable.
\end{corollary}

It was proved in \cite[Proposition 5.1]{ZH2020} that a pseudo-$\aleph_{1}$-compact and $\mathcal{M}$-factorizable topological group is $\mathbb{R}$-factorizable. Moreover, it is well-known that a topological group is $m$-factorizable if and only if it is $\mathbb{R}$-factorizable and pseudo-$\aleph_{1}$-compact.

\begin{proposition}\label{6mt999}
Every pseudo-$\aleph_{1}$-compact $\mathcal{M}$-factorizable topological group is $S\mathcal{M}$-factorizable.
\end{proposition}

During the following figure, for convenience, $\mathbb{R}$ denotes $\mathbb{R}$-factorizable topological groups, $\mathcal{M}$ denotes $\mathcal{M}$-factorizable topological groups, and so on. Moreover, the two arrows represent stronger properties, such as every $P\mathbb{R}$-factorizable topological group is $\mathbb{R}$-factorizable, .etc. The notions $+nar$, $+fea$ and $+pse$ represent adding the properties of $\omega$-narrow, feathered and pseudo-$\aleph_{1}$-compact, respectively, and then the single arrow will hold, that is, two spaces at both ends of the arrow will be equivalent.

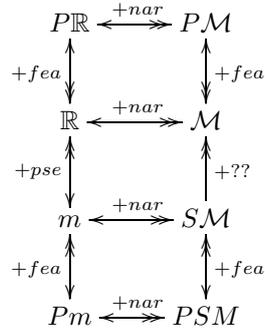
\begin{equation*}
\centerline{
\xymatrix{ P\mathbb{R} \ar@{>>}[d]_{+fea} \ar@{>>}[r]^{+nar} &P\mathcal{M} \ar[l] \ar@{>>}[d]^{+fea} \\
  \mathbb{R} \ar[u] \ar[d] \ar@{>>}[r]^{+nar} & \mathcal{M} \ar[l] \ar[u]   \\
  m  \ar[d] \ar@{>>}[u]^{+pse} \ar@{>>}[r]^{+nar} & S\mathcal{M} \ar[d] \ar[l] \ar@{>>}[u]_{+??}    \\
  Pm \ar@{>>}[u]^{+fea} \ar@{>>}[r]^{+nar} & PSM \ar[l] \ar@{>>}[u]_{+fea}        }}
\end{equation*}

However, the property of pseudo-$\aleph_{1}$-compact is too strong in Proposition \ref{6mt999}, so it is natural to consider what properties can be added to an $\mathcal{M}$-factorizable topological group to imply it is $S\mathcal{M}$-factorizable and every $S\mathcal{M}$-factorizable topological group obtains the properties at the same time.

Finally, the following question is posed naturally.

\begin{question}
Is there an $S\mathcal{M}$-factorizable topological group but not $m$-factorizabe?
\end{question}

\end{document}